\documentclass[11pt]{amsart}

\usepackage{amsmath,amssymb,amscd,latexsym}
\usepackage[sans]{dsfont}

\title{Generators of simple Lie algebras II}
\author{Jean-Marie Bois}
\date{26 June 2008}


\address{Fakult\"at f\"ur Mathematik, Universit\"at Bielefeld, Postfach 10 01 31, 33501 Bielefeld, Germany}
\email{jbois@math.uni-bielefeld.de}

\newenvironment{thm}{\noindent \bf Theorem.\it}{ \rm }
\newenvironment{prop}{\noindent \bf Proposition.\it}{ \rm }
\newenvironment{lemma}{\noindent \bf Lemma.\it}{ \rm }
\newenvironment{cor}{\noindent \bf Corollary.\it}{ \rm }

\renewcommand{\proof}{\noindent {\bf Proof.} }
\newcommand{\nr}{n$^o$}
\newcommand{\nbf}{\noindent \bf }

\newcommand{\FF}{\mathbb{F}}
\newcommand{\NN}{\mathbb{N}}
\newcommand{\ZZ}{\mathbb{Z}}

\newcommand{\gtm}{\mathfrak{m}}

\newcommand{\cO}{\mathcal{O}}
\newcommand{\cN}{\mathcal{N}}
\newcommand{\cS}{\mathcal{S}}

\newcommand{\nn}{\boldsymbol{\underline{n}}}
\newcommand{\un}{\mathds{1}}
\newcommand{\id}{\operatorname{id}}

\newcommand{\Der}{\operatorname{Der}}
\newcommand{\End}{\operatorname{End}}
\newcommand{\ad}{\operatorname{ad}}

\renewcommand{\d}{\partial}
\renewcommand{\epsilon}{\varepsilon}
\renewcommand{\div}{\operatorname{div}}

\newcommand{\Fxy}[1]{\FF \langle #1 \rangle}

\begin{document}

\renewcommand{\labelitemi}{\textbullet}
\renewcommand{\labelenumi}{\arabic{enumi}.}

\makeatletter \def\revddots{\mathinner{\mkern1mu\raise\p@ \vbox{\kern7\p@\hbox{.}}\mkern2mu \raise4\p@\hbox{.}\mkern2mu\raise7\p@\hbox{.}\mkern1mu}} \makeatother 

\maketitle

\begin{abstract}
This paper is a continuation of earlier work on generators of simple Lie algebras in arbitrary characteristic \cite{jmb}. We show that, in contrast to classical Lie algebras, simple graded Lie algebras of Cartan type $S,H$ or $K$ never enjoy the ``one-and-a-half generation'' property. The methods rely on a study of centralisers in Cartan type Lie algebras.
\end{abstract}

\section*{Introduction}

In this paper we continue the study of generators of simple Lie algebras initiated in \cite{jmb}. The original motivation for this question comes from group theory. Namely, any finite simple group can be generated by 2 elements; furthermore, one of these elements can be chosen arbitrarily (see for example \cite{Aschbacher,Guralnick}). This last property is referred to as the ``one-and-a-half generation property''. Recently, B. Kunyavski\u\i\ asked whether the analogous property holds for simple Lie algebras (quoted in \cite[Problem 5]{Premet}).

In \cite{jmb}, we provided a partial answer to these questions. We first proved that any simple Lie algebra over an algebraically closed field $\FF$ of characteristic $\neq 2,3$ can be generated by 2 elements. For the 1,5-generation property we have the following information. First, any classical simple Lie algebra (ie. simple complex Lie algebras and their modular analogues) is generated by 1,5 elements. Also, the graded Cartan type Lie algebra $W(m,\nn)$ has the 1,5 generation property if and only if $m = 1$. In the present work we show that the simple graded Cartan type Lie algebras of the remaining types $S,H$ or $K$ are never generated by 1,5 elements. 

The methods used in \cite{jmb} to investigate the 1,5 generation property were geometrical in the classical case, and {\it ad hoc} for algebras of Cartan type $W$. In the present paper the approach involves properties of centralisers in filtered Lie algebras. It should be noted that the criterion stated below does not apply for the Lie algebras $W(m,\nn)$.

This paper is organised as follows. We first give a criterion for the 1,5 generation property to fail to be satisfied in a Lie algebra. In the special case of filtered Lie algebras we obtain the following: \\

{\nbf Criterion.} \it 
Let $L = L_{\geq -r} \supseteq \ldots \supseteq L_{\geq s} \supseteq (0)$ be a decreasingly filtered Lie algebra, with $r,s \geq 0$. Let $H \subseteq L$ be an ideal of $L$ such that the centraliser $C_L(H) = (0)$. Assume that, for every $h \in H$, the centraliser of $h$ in $L$ intersects the positive part $L_{\geq 1}$ non-trivially. Then $H$ cannot be generated by 1,5 elements. \\
\rm

We proceed with a brief analysis of centralisers in the Lie algebra $W_m = \Der(B_m)$ of derivations of a truncated polynomial ring in $m$ variables. We give a normal form for derivations $D \in W_m$ such that the ring of $D$-constants is $B_m^D = \FF$. Using this normal form we obtain that for such derivations, the centraliser $W_m^D$ has dimension $m$, the lowest possible value.

In Section \ref{section2} we deal with the case of simple graded Lie algebras of Cartan type $S,H$ and $K$. In each case we prove failure of the 1,5-generation property by application of the above criterion on centralisers. \\

Throughout, we will work over an algebraically closed field $\FF$ of characteristic $p > 3$. All maps, derivations, etc.\ under consideration will be $\FF$-linear. We will omit the corresponding subscripts and write for example $\Der$ or $\otimes$ instead of $\Der_\FF$ or $\otimes_\FF$.

\section{Preliminaries}

\subsection{A criterion for one-and-a-half generation}

\subsubsection{} \label{criterion}
Let $L$ be a Lie algebra. Recall that $L$ is {\it generated by 1,5 elements} if, for all $x \in L$ with $x \neq 0$, there exist $y \in L$ such that the generated subalgebra $\Fxy{x,y} = L$. We have the following criterion for the 1,5-generation property to fail: \\

\begin{prop} Let $L$ be a Lie algebra. The following are equivalent:

(i) $L$ is {\bf not} generated by 1,5 elements;

(ii) there exist a dense subset $\Omega \subseteq L$ and a family of subalgebras $(L_i)_{i \in I}$ such that:
\begin{itemize}
\item for all $i \in I$, we have $L_i \neq L$;
\item the union $\bigcup_{i \in I} L_i \supseteq \Omega$;
\item the intersection $\bigcap_{i \in I} L_i \neq (0)$. \\
\end{itemize}
\end{prop}

\proof\ Assume first that $L$ is not generated by 1,5 elements: there exists a nonzero $x \in L$ such that, for all $y \in L$, the subalgebra $L_y = \Fxy{x,y} \neq L$. We have $0 \neq x \in \bigcap_{y \in L} L_y$. Also, $\Omega = L = \bigcup_{y \in L} \{ y \} \subseteq \bigcup_{y \in L} L_y$. Thus property $(ii)$ holds.

Now assume that $(ii)$ is satisfied, and let us show that $(i)$ holds. Let $0 \neq x \in \bigcap_{i \in I} L_i$. Let $U = \{ y \in L\ |\ \Fxy{x,y} = L\}$. This is an open subset of $L$. If $L$ was generated by 1,5 elements, then we would have $U \neq \emptyset$. Since $\Omega$ is dense, there exists an element $y \in \Omega \cap U$. Because $y \in U$, we have $\Fxy{x,y} = L$. On the other hand, since $y \in \Omega$, there exists $i \in I$ such that $y \in L_i$. Since also $x \in L_i$, we have $\Fxy{x,y} \subseteq L_i \neq L$, a contradiction.

\subsubsection{} \label{corollary1,5}
For a Lie algebra $L$, we will use the notation $L^X$ for the centraliser of $X$ in $L$, where $X$ can be an element or a subspace of $L$. \\

\begin{cor}
Let $L = L_{\geq -r} \supseteq \ldots \supseteq L_{\geq s} \supseteq (0)$ be a decreasingly filtered Lie algebra, with $r,s \geq 0$. Let $H \subseteq L$ be an ideal of $L$ satisfying the following conditions:
\begin{itemize}
\item the centraliser $L^H = (0)$;
\item there exists a dense subset $\Omega \subseteq H$ such that, for all $X \in \Omega$, $L^X \cap L_{\geq 1} \neq (0)$.
\end{itemize}
Then $H$ is not generated by 1,5 elements. \\
\end{cor}

{\nbf Remark.} Note that the second condition concerns the centraliser of $X$ in the full algebra $L$, not only in $H$. \\

\proof\ For all $X \in \Omega$, we choose a non-zero element $Y \in L_{\geq 1}$ such that $[X,Y] = 0$. Denote by $H_X = L^Y \cap H$, the space of elements in $H$ which commute to $Y$. Then $X \in H_X$ by construction, so that $\bigcup_{X \in \Omega} H_X \supseteq \Omega$. Since $Y \neq 0$ centralises $H_X$, we must have $H_X \neq H$ according to the first condition. Now, the subalgebra $H$ inherits a decreasing filtration from the one in $L$; let $H_{\geq s'}$ be the non-zero subspace of maximal index. Since $Y \in L_{\geq 1}$, we have $[Y,H_{\geq s'}] \subseteq H_{> s'} = (0)$, whence $\bigcap_{X \in \Omega} H_X \supseteq H_{\geq s'} \neq (0)$. The corollary now follows from the proposition above.

\subsubsection{\bf Remark.} For the Jacobson-Witt algebra $L = H = W_m$ (and, more generally, for the non-restricted $W(m,\nn)$) it can be shown that the second condition of Corollary \ref{corollary1,5} does not hold. Indeed, let $U \subseteq L$ be the defined as:
\begin{equation*}
U = \{ X \in L \ | \ L^X \cap L_{\geq 0} = (0) \}.
\end{equation*}
Equivalently, $X \in U$ if and only if the linear map $(\ad\, X)_{| L_{\geq 0}}: L_{\geq 0} \to L$ is injective. From this it is easy to see that $U$ is open. Let $\cN$ be the regular nilpotent derivation as in \cite[Lemma 3]{Premet_Wn}: using \cite[Lemma 7 (ii)]{Premet_Wn} we can see that $\cN \in U$, so that $U$ is not empty. This precludes the existence of a dense subset $\Omega \subseteq L$ satisfying the second property of the previous corollary.

\subsection{Derivations of truncated polynomial rings}

\subsubsection{} \label{TPrings}
Let $m \geq 1$ be an integer. Let $B_m$ denote the truncated polynomial ring in $m$ variables over $\FF$, that is:
\begin{equation}
B_m = \FF[x_1,\ldots,x_m] / (x_1^p,\ldots,x_m^p).
\end{equation}
The generators of $B_m$ are again denoted by $x_1,\ldots,x_m$, so that the equality $x_i^p = 0$ holds in $B_m$. Recall that $\dim(B_m) = p^m$.

Let $W_m = \Der(B_m)$ be the full derivation algebra. For any $D \in W_m$, we can consider the subring of $D$-constants $B_m^D \subseteq B_m$ to be the space of truncated polynomials annihilated by $D$. We always have $\FF \subseteq B_m^D$.

\subsubsection{}
We introduce some convenient terminology. A derivation $D \in W_m$ is {\it regular nilpotent} if $D^{p^m-1} \neq D^{p^m} = 0$ and {\it regular semi-simple} if $D$ has $p^m$ distinct eigenvalues as a linear map acting on $B_m$. \\

\begin{prop}
Let $D \in W_m$ be a derivation of $B_m$. Assume that $D$ has no non-trivial constants, ie. $B_m^D = \FF$. Then there exists a decomposition $B_m = B' \otimes B''$, with $B' \simeq B_r$ and $B'' \simeq B_{m-r}$, and derivations $\cN \in \Der (B') \simeq W_r$, $\cS \in \Der (B'') \simeq W_{m-r}$ satisfying the following properties:
\begin{enumerate}
\item $\cN \in W_r$ is regular nilpotent and $\cS \in W_{m-r}$ is regular semi-simple;
\item $D = \cN \otimes \id_{B''} + \id_{B'} \otimes \cS \in \Der(B' \otimes B'')$. \\
\end{enumerate}
\end{prop}

\proof\ Let $\chi(t)$ be the characteristic polynomial of $D$ acting in $B_m$ and $\Lambda$ be the set of roots of $\chi(t)$. By \cite[Theorem 1, (ii)]{Premet_Wn}, $\chi(t)$ is a $p$-polynomial, and therefore all its roots have same multiplicity which is a power of $p$, say $p^r$. Furthermore, the set $\Lambda$ is an additive subgroup of $\FF$, of order $p^{m-r}$.

Let $B'$ be the nilspace of $D$, that is the kernel of $D^{p^r}$, and $\cN$ be the derivation of $B'$ induced by $D$. Then $\cN$ is nilpotent and has one-dimensional kernel $\FF$. Hence, there exists a linear basis $\{ z_1,\ldots , z_{p^r}\}$ of $B'$ such that $\cN(z_k) = z_{k-1}$ for all $k > 1$ and $\cN(z_1) = 0$; we may assume that $z_1 = 1$. It is easy to see that the algebra $B'$ has no non-trivial ideal which is stable under $\cN$: a fortiori no non-trivial ideal is stable under all derivations of $B'$. By Harper's theorem \cite{Harper}, $B'$ is isomorphic to some truncated polynomial ring; because $\dim B' = p^r$ we must have $B' \simeq B_r$.

Let $\lambda \in \Lambda$ be an eigenvalue of $D$. Denote by $B^\lambda$ be the corresponding eigenspace and $B^{(\lambda)}$ the characteristic subspace; in particular $\dim(B^{(\lambda)}) = p^r$. Also, let $g_{\lambda}$ be a non-zero eigenvector; for $\lambda = 0$ we can choose $g_0 = 1$. Finally, for all $k \in \{1,\ldots,p^r\}$, set $g_{\lambda,k} = z_k g_{\lambda}$. We have $D(g_{\lambda,k}) = \lambda\, g_{\lambda,k} + g_{\lambda,k-1}$. From this fact we can see that the matrix of the action of $D$ on $B^{(\lambda)}$ is given by a single Jordan block of size $p^r$:
\begin{equation} 
\left[ \begin{array}{cccc}
\lambda & 1 & \ldots & 0 \\
0 & \ddots & \ddots  & \\
0 &  & \lambda & 1 \\
0 &  \ldots & 0 & \lambda 
\end{array} \right] 
\end{equation}
and the elements $g_{\lambda,k}$ for $k \in \{ 1,\ldots,p^r\}$ span the subspace $B^{(\lambda)}$. Furthermore, we have $\dim ( B^{\lambda} ) = 1$.

Let $B'' = \sum_{\lambda \in \Lambda} B^{\lambda}$: clearly, $B''$ is a subalgebra of $B_m$. Since the $g_{\lambda,k} = z_k g_\lambda$, where $\lambda \in \Lambda$ and $k \in \{1,\ldots,p^r \}$, form a linear basis of $B_m$, we can see that the multiplication map induces an isomorphism $B' \otimes B'' \stackrel{\sim}{\to} B_m$. Since $B''$ is generated by eigenvectors of $D$, it is also stable under the derivation $D$. Let $\cS \in \Der(B'')$ be the induced derivation. By construction, $\cS$ is semi-simple, and all eigenspaces of $\cS$ in $B''$ have dimension 1, so that $\cS$ is also regular. Finally, it is easy to check that $D = \cN \otimes \id_{B''} + \id_{B'} \otimes \cS$.

To complete the proof of the proposition it remains to show that $B'' \simeq B_{m-r}$. Let $\gtm \subseteq B_m$ be the maximal ideal of $B_m$. Let $\xi_1,\ldots,\xi_r \in B'$ be a set of generators of $B'$ given by an isomorphism $B_r = \FF[x_1,\ldots,x_r] / (x_1^p,\ldots,x_r^p) \stackrel{\sim}{\to} B'$. Using the fact that $(\xi_1 \cdots \xi_r)^{p-1} \neq 0$, we can see that the linear terms of $\xi_1,\ldots,\xi_r$ are linearly independent modulo $\gtm^2$. Hence, they admit a cobasis $\{ \xi_{r+1},\ldots,\xi_m \}$ for $\gtm$ modulo $\gtm^2$. We obtain a chain of isomorphisms:
\begin{equation}
B_{m-r} \simeq \FF[\xi_{r+1},\ldots,\xi_m] \simeq \frac{B_m}{ \sum_{i = 1}^r \xi_i\, B_m } \simeq \left( \frac{B'}{ \sum_{i = 1}^r \xi_i\, B' } \right) \otimes B'' \simeq \FF \otimes B'' \simeq B''.
\end{equation}
The proposition is proved.

\subsubsection{} \label{CorSmallKernel}
\begin{cor}
Let $D \in W_m$ be such that $B_m^D = \FF$ and $W_m^D$ be the centraliser of $D$ in $W_m$. Then $\dim ( W_m^D ) = m$. \\
\end{cor}

\proof\ We first prove the corollary in two particular cases. First, if $D = \cN$ is regular nilpotent, the equality $\dim ( W_m^D ) = m$ follows from \cite[Lemma 7.2 (ii)]{Premet_Wn}.
Second, assume that $D = \cS$ is regular semi-simple. Then $D$ is contained in some maximal torus $T \subseteq W_m$. Since $D$ is regular, we can see that $T = \sum_{j \geq 0}\FF\, D^{[p^j]}$, and $W_m^D = W_m^T$ the centraliser of $T$. By \cite[Corollary 7.5.17]{Strade}, $W_m^T$ is a Cartan subalgebra of minimal dimension, so that $\dim (W_m^T) = m$ by \cite[Table 3.25]{Benkart}.

Now we turn to the general case. We write $B_m = B' \otimes B''$ and $D = \cN \otimes \id_{B''} + \id_{B'} \otimes \cS$ as in the Proposition. Let $\Delta \in W_m^D$. From the construction in the proof, we have that $B'$ is the nilspace of $D$ and $B''$ the sum of eigenspaces, so that $\Delta$ must stabilise both $B'$ and $B''$. This means that $\Delta$ takes the form $\Delta = \Delta' \otimes \id_{B''} + \id_{B'} \otimes \Delta''$, where $\Delta' \in \Der(B') \simeq W_r$ and $\Delta'' \in \Der(B'') \simeq W_{m-r}$. Furthermore, since $[D,\Delta] = 0$ we also have $[\Delta', \cN] = 0$ and $[\Delta'',\cS] = 0$. It follows readily that:
\begin{equation}
W_m^D = W_r^{\cN} \otimes \id_{B''} +  \id_{B'} \otimes W_{m-r}^{\cS}.
\end{equation}
Since the corollary holds in the regular nipotent and semi-simple cases, we obtain $\dim (W_m^D) = \dim (W_r^{\cN}) + \dim ( W_{m-r}^{\cS} ) = r + (m - r) = m$ as we wanted.

\section{Generators for graded Cartan type algebras}
\label{section2}

\subsection{Divided power algebras and special derivations} $ $\\ 
\label{SubSectionW}

We recall from \cite[Section 2.1]{Strade} some basic facts and notations about special derivations of divided power algebras.

\subsubsection{}
Let $m \geq 1$ be an integer. For a multi-index $\alpha=(\alpha_1,\ldots,\alpha_m) \in \ZZ^m$, we define the length to be $|\alpha| = \alpha_1+\ldots+\alpha_m$. The factorial and binomial coefficients are given, for $\alpha,\beta \in \NN^m$, by:
\begin{equation}
\alpha! = \alpha_1 ! \cdots \alpha_m ! \quad \mbox{ and } \quad \binom{\alpha}{\beta} = \binom{\alpha_1}{\beta_1} \cdots \binom{\alpha_m}{\beta_m}.
\end{equation}
For all $k \in \{1,\ldots,m\}$, define the $k$-th unit multi-index by $\epsilon_k = (\delta_{1,k},\ldots,\delta_{m,k})$. Furthermore, we set $\un = \epsilon_1 + \ldots + \epsilon_m = (1,\ldots,1)$. The set of multi-indices can be partially ordered in the following fashion: for $\alpha,\beta \in \ZZ^m$, let $\alpha \leq \beta$ if the components $\alpha_i \leq \beta_i$ for all $i \in \{1,\ldots,m \}$. We will also write $\alpha < \beta$ if $\alpha \leq \beta$ and $\alpha \neq \beta$.

\subsubsection{\bf Conventions} 
\label{furtherconventions}
In the sequel we will use the following notational conventions. For a multi-index denoted by $\nn \in \NN^m$ we always make the implicit assumption that all components $n_i > 0$. Furthermore, given such a multi-index we also implicitly define $\tau = \tau(\nn) \in \NN^m$ to be $\tau = (p^{n_1}-1,\ldots,p^{n_m}-1)$. In particular, the expression $0 \leq \alpha \leq \tau$ means that the components $0 \leq \alpha_i \leq p^{n_i} -1$ for all $i \in \{ 1,\ldots,m\}$.

\subsubsection{\bf Divided power algebras} 
\label{dividedpowers}
Let $m \geq 1$ be an integer. Define the divided power algebra $\cO(m)$ over $\FF$ as follows: as vector space, $\cO(m) = \FF[x^{(\alpha)}\, :\, \alpha \in \NN^m]$ and the multiplication map is given by the formula:
\begin{equation}
(\forall\, \alpha,\beta \in \NN^m)\ :\ x^{(\alpha)} x^{(\beta)} = \binom{\alpha + \beta}{\alpha} x^{(\alpha+\beta)}.
\end{equation}
The algebra $\cO(m)$ has a natural grading over $\ZZ$ such that all elements $x^{(\alpha)}$ are homogeneous, of degree $| \alpha |$. For all muti-indices $\nn \in \NN^m$ , define a linear subspace:
\begin{equation}
\cO(m,\nn) = \operatorname{Span} \left\{x^{(\alpha)}\ |\ 0 \leq \alpha\leq \tau \right\} \subseteq \cO(m).
\end{equation}
Then $\cO(m,\nn)$ is actually a graded subalgebra of $\cO(m)$, of dimension $p^{|\nn|}$. 

\subsubsection{} \label{Suggestive}
For divided power algebras, we will use suggestive terminology by analogy with polynomial rings. For example, elements $x^{(\alpha)}$ will be called ``monomials in $x_1,\ldots,x_m$, and a monomial $x^{(\alpha)}$ is ``independent of $x_k$'' if the component $\alpha_k = 0$.

We will also use the following notations. For all $k \in \{ 1,\ldots,m\}$ and $r \geq 0$, set $x_k^{(r)} = x^{(r\, \epsilon_k)}$; when $r = 1$, write $x_k = x_k^{(1)}$. With this notation, we can check that $x^{(\alpha)} = x_1^{(\alpha_1)} \cdots x_m^{(\alpha_m)}$, and $x_k^{(r)} x_k^{(s)} = \binom{r+s}{r} x_k^{(r+s)}$.

\subsubsection{\bf Special derivations}  
For all $k \in \{1,\ldots,m\}$, we define the ``partial derivative'' $\d_k \in \Der\, \cO(m,\nn)$ by the rule $\d_k(x^{(\alpha)}) = x^{(\alpha - \epsilon_k)}$. The Lie algebra of special derivations of $\cO(m,\nn)$, denoted by $W(m,\nn)$, is defined to be the sub-$\cO(m,\nn)$-module of $\Der\, \cO(m,\nn)$ generated by these partial derivatives:
\begin{equation} \label{plushaut}
W(m,\nn) = \cO(m,\nn)\, \d_1 \oplus \ldots \oplus  \cO(m,\nn)\, \d_m \subseteq \Der\, \cO(m,\nn).
\end{equation}
The inclusion is proper unless $\nn = \un$. In this case, we have $\cO(m,\un) \simeq B_m$ and $W(m,\un) \simeq W_m$, where $B_m$ and $W_m$ are defined in Section \ref{TPrings}.

\subsubsection{} 
Because $\cO(m,\nn)$ is a graded algebra, the Lie algebra $W(m,\nn)$ inherits a natural grading, compatible with the grading on $\cO(m,\nn)$. More precisely, all derivations of the form $x^{(\alpha)} \d_k$ have degree $|\alpha| - 1$. The homogeneous subspace of maximal degree has dimension $m$ and is generated by derivations $x^{(\tau)}\, \d_k$, for $\tau = (p^{n_1}-1,\ldots,p^{n_m}-1)$ and $k \in \{1,\ldots,m\}$. The homogeneous subspace of minimal degree -1 has also dimension $m$ and is generated by the partial derivatives $\d_k$, for $k \in \{1,\ldots,m\}$.

We will frequently use natural notations such as $W(m,\nn)_{\geq k}$ to denote the subspace spanned by homogeneous elements of degree $\geq k$, or $W(m,\nn)_{max}$ to denote the homogeneous subspace of maximal degree. Also, we will use expressions such as ``$D = D_r + \ldots$, where ``$\ldots$'' are higher order terms''. This will mean that $D_r$ is homogeneous of degree $r$, and ``$\ldots$'' is an element of $W(m,\nn)_{> r}$.

\subsubsection{\bf Embeddings of $W$ type algebras.}
\label{Plongements}
Let $m,\nn$ be as above and let $N = |\nn| \in \NN$. By \cite[p. 64]{Strade}, the algebra $\cO(m,\nn)$ is isomorphic to $\cO(N,\un)$ as abstract (ie. non-graded) algebras. Any isomorphism $\sigma:\cO(m,\nn) \stackrel{\sim}{\to} \cO(N,\un)$ gives rise to an isomorphism $\Der\, \cO(m,\nn) \stackrel{\sim}{\to} \Der\, \cO(N,\un) = W(N,\un)$, hence an embedding:
\begin{equation}
\iota: W(m,\nn) \hookrightarrow W(|\nn|,\un).
\end{equation}
Note that, like $\sigma$, this embedding does not respect the natural gradings. By construction, we have the following relations, for all $D \in\, W(m,\nn)$ and $f \in\, \cO(m,\nn)$:
\begin{equation}
\iota(D) = \sigma^{-1} \circ D \circ \sigma
\end{equation}
\begin{equation} \label{iotapreservesmodulestructure}
\iota(f \cdot D) = \sigma(f) \cdot \iota(D).
\end{equation}

\subsection{Lie algebras of type $S$.} $ $\\ 
\label{SubSectionS}

In this section we focus on simple graded Lie algebras of type $S$. We first recall some constructions; the question of generators is addressed in Theorem \ref{UnEtDemiS}.

\subsubsection{}
Let $m \geq 2$ be an integer and $\nn$ a multi-index. Consider the Lie algebra $W(m,\nn)$ as in \ref{SubSectionW}. Define the {\it divergence of a special derivation} by the rule:
\begin{equation}
\div (\sum_{j = 1}^m f_j\, \d_j) = \sum_{j=1}^m\, \d_j(f_j) \in \cO(m,\nn).
\end{equation}
The divergence map satisfies the following properties, for all $f \in\, \cO(m,\nn)$ and $D,D' \in\, W(m,\nn)$:
\begin{equation} \label{divfD}
\div(f\, D) = f\, \div(D) + D(f),
\end{equation}
\begin{equation} \label{divDD'}
\div\, [D,D'] = D(\div\, D') - D'(\div\, D).
\end{equation}
The Lie algebra $S(m,\nn) \subseteq W(m,\nn)$ is defined as the subspace of divergence-free derivations, ie. $D \in S(m,\nn)$ if and only if $\div (D) = 0$. Then, the derived subalgebra $S(m,\nn)^{(1)}$ is a simple ideal of codimension $m$. Furthermore, both algebras $S(m,\nn)$ and $S(m,\nn)^{(1)}$ are graded subalgebras of $W(m,\nn)$.

For simple Lie algebras of type $S$, it is customary to impose the further restriction $m \geq 3$; for $m = 2$ the Lie algebra $S(2,\nn)^{(1)}$ can actually be identified with a simple Lie algebra of type $H$. However, the results we will prove here do not require such a restriction, so that the case $m = 2$ will also be considered. 

\subsubsection{} \label{IotaRespectsDiv}
Recall the embedding $\iota: W(m,\nn) \hookrightarrow W(N, \un)$ described in Section \ref{Plongements}), with $N = \nn$. The following lemma ensures that we can arrange $\iota$ to preserve the divergence maps. \\

\begin{lemma}
There exists an isomorphism $\sigma: \cO(m,\nn) \stackrel{\sim}{\to} \cO(N,\un)$ satisfying the following property. Let  $\iota: W(m,\nn) \hookrightarrow W(N, \un)$ the induced embedding (see \ref{Plongements}). Denote by $\div_{(m,\nn)}$ and $\div_{(N,\un)}$ the divergence maps defined in $W(m,\nn)$ and $W(N,\un)$ respectively. Then, for all $D \in W(m,\nn)$:
\begin{equation} \label{CompatibilityIotaDiv}
\div_{(N,\un)} \big( \iota(D) \big) = \div_{(m,\nn)}(D).
\end{equation}
\end{lemma}

\proof\ 
We use the usual notations for elements of $\cO(m,\nn)$ and $W(m,\nn)$. For $\cO(N,\un)$ and $W(N,\un)$ we introduce alternative notations. To make descriptions easier, we won't use the standard variables indexed by $\{1,\ldots,N\}$ but rather some suitably doubly indexed variables. We will consider $\cO(N,\un)$ to be generated by variables $\xi_{i,j} = \xi_{i,j}^{(1)}$, where the pairs $(i,j)$ satisfy $i \in \{1,\ldots,m\}$ and $j \in \{ 0,\ldots,n_i - 1\}$. Any element $\Delta \in W(N,\un)$ can be written as:
\begin{equation}
\Delta = \sum_{i,j}\, \varphi_{i,j}\, \delta_{i,j},
\end{equation}
where the coefficients $\varphi_{i,j} \in \cO(N,\un)$ and $\delta_{i,j}$ is the partial derivation defined by $\delta_{i,j}(\xi_{i,j}) = 1$ and $\delta_{i,j}(\xi_{k,l}) = 0$ otherwise. For such a derivation, we have:
\begin{equation}
\div_{(N,\un)}(\Delta) = \sum_{i,j}\, \delta_{i,j}(\varphi_{i,j}).
\end{equation}
The isomorphism $\sigma: \cO(m,\nn) \stackrel{\sim}{\to} \cO(N,\un)$ is uniquely defined by the assignment $x_i^{(p^j)} \mapsto \xi_{i,j}$.

By construction, the embedding $\iota$ preserves multiplication by ``polynomials'', see Equation  (\ref{iotapreservesmodulestructure}). Using the formula  (\ref{divfD}) concerning $\div(fD)$, it is clearly enough for the Lemma to prove the identity (\ref{CompatibilityIotaDiv}) on a generating system of the $\cO(m,\nn)$-module $W(m,\nn)$, for example the partial derivatives $\d_1,\ldots,\d_m$.

Write out the derivation $\iota(\d_k) = \sum_{i,j}\, \psi_{i,j} \, \delta_{i,j}$. In this expression, we have $\psi_{i,j} = \iota(\d_k)(\xi_{i,j}) = \sigma\big( \d_k(x_i^{(p^j)}) \big)$, so that:
\begin{equation}
\psi_{i,j} = 
\left\{ \begin{array}{cl}
0 & \mbox{ if }  i \neq k, \\
\sigma(x_k^{(p^j-1)}) & \mbox{ if } i = k.
\end{array} \right.
\end{equation}
Using the decomposition $p^j-1 = (p-1) + (p-1)p + \ldots + (p-1)p^{j-1}$, we can see that 
\begin{equation}
x_k^{(p^j-1)} \in \FF\, \left(x_k^{(1)}\right)^{p-1} \cdots \left( x_k^{(p^{j-1})}\right)^{p-1},
\end{equation}
so that $\psi_{k,j} = \sigma(x_k^{(p^j-1)}) \in \FF\, \xi_{k,0}^{p-1} \cdots \xi_{k,j-1}^{p-1}$. It readily follows that $\delta_{k,j} (\psi_{k,j}) = 0$; we finally obtain:
\begin{equation}
\div_{(N,\un)}(\iota(\d_k)) = \sum_{i,j}\, \, \delta_{i,j}(\psi_{i,j}) = 0 = \div_{(m,\nn)}(\d_k),
\end{equation}
as desired. The lemma is proved.

\subsubsection{} \label{UnEtDemiS}
\begin{thm}
Let $m \geq 2$ be an integer and $\nn \in \NN^m$ a multi-index with non-zero entries. The simple Lie algebra $S(m,\nn)^{(1)}$ is not generated by 1,5 elements. \\
\end{thm}

\proof\ First we show that, for all $D \in S(m,\nn)^{(1)}$, the constant ring $\cO(m,\nn)^D \neq \FF$. By Lemma \ref{IotaRespectsDiv}, we have an embedding $S(m,\nn)^{(1)} \hookrightarrow S(N,\un)^{(1)}$ induced by some isomorphism $\cO(m,\nn) \simeq \cO(N,\un)$, so that we may assume that $\nn = \un$. In this case we use the simpler notations $W_m$, $S_m$, $B_m$ instead of $W(m,\un)$, $S(m,\un)$, $\cO(m,\un)$. Let $CS_m = \{ D \in W_m |\ \div(D) \in \FF \}$. We can see that $S_m^{(1)}$ is an ideal of codimension $m+1$ in $CS_m$. Thus, for every $D \in S_m^{(1)}$ the subspace $[D,CS_m] \subseteq CS_m$ has codimension $\geq m+1$; therefore, the centraliser $CS_m^D$ has dimension $\geq m+1$. A fortiori, we have $\dim (W_m^D) > m$: by Corollary \ref{CorSmallKernel}, we obtain $\dim\, B_m^D > 1$.
 
Now we go back to the general case of $S(m,\nn)^{(1)}$. We will apply Corollary \ref{corollary1,5} with $L = S(m,\nn)$ and $H = S(m,\nn)^{(1)}$. The fact that the centraliser of $H$ in $L$ is trivial follows from the fact that both $H$ and $L$ have trivial centre. Now we turn to the second condition. We endow $S(m,\nn)$ with the filtration defined by its standard grading. Let $\Omega \subseteq S(m,\nn)^{(1)}$ be the set of derivations which have a non-zero homogeneous part of degree $-1$; this is a non-empty Zariski open subset of $S(m,\nn)^{(1)}$. All we need to show is that, for all $D \in \Omega$, there exists a non-zero element $\Delta \in S(m,\nn)_{\geq 1}$ such that $[D,\Delta] = 0$ (we do not require $\Delta$ to belong to the derived subalgebra).

We let $D \in \Omega$ act as a derivation on $\cO(m,\nn)$. Let us check that $\cO(m,\nn)_{\geq 2}^D \neq (0)$, ie. there exist a non-zero element $f \in \cO(m,\nn)$ such that $D(f) = 0$ and $f$ has no homogeneous terms of degree 0 and 1. We showed above that $\cO(m,\nn)^D \neq \FF$, so that there exist $0 \neq f \in \cO(m,\nn)^D$ with no constant term. Write out $f = f_1 + \ldots$, where $f_1$ is homogeneous of degree 1 and ``$\ldots$'' are higher order terms. If $f_1 = 0$, we have $f \in \cO(m,\nn)_{\geq 2}$ and we are done. If $f_1 \neq 0$, we have $f^2 = f_1^2 + \ldots \neq 0$, and we can use $f^2 \in \cO(m,\nn)_{\geq 2}^D$.

Now choose a non-zero element $f \in \cO(m,\nn)_{\geq 2}^D$, and let $\Delta = f\, D \in W(m,\nn)$. Because $D$ has a non-zero term of degree $-1$ and $f \in \cO(m,\nn)_{\geq 2}$, we have $0 \neq \Delta = f\, D \in W(m,\nn)_{\geq 1}$. Using the formula (\ref{divfD}) we obtain $\div(\Delta) = 0$, hence $\Delta \in S(m,\nn)_{\geq 1}$. And finally, we have $[D,\Delta] = [D,f\, D] = 0$. This is what remained to be shown: thus, we can apply Corollary \ref{corollary1,5} and obtain that $S(m,\nn)^{(1)}$ is not generated by 1,5 elements.

\subsection{Lie algebras of type $H$.} $ $\\ \label{SubSectionH}

In this section, we let $m = 2r \geq 2$ be an even integer and $\nn \in \NN^{2r}$ a multi-index with non-zero entries. We keep the general notations for divided power algebra and special derivations introduced in Section \ref{SubSectionW}.

\subsubsection{} \label{introH}
We recall some construction from \cite[Section 4.2]{Strade}. First define a binary operation on $\cO(2r,\nn)$ by setting:
\begin{equation}
\{ f,g \} = \sum_{j=1}^{2r} \sigma(j) \d_j(f)\, \d_{j'}(g),
\end{equation}
where $j' = \left\{ \begin{array}{rcl} j+r & \mbox{if} & 1 \leq j \leq r \\ j-r & \mbox{if} & r < j \leq 2 r \end{array} \right.$ and $\sigma(j) = \left\{ \begin{array}{rcl} +1 & \mbox{if} & 1 \leq j \leq r \\ -1 & \mbox{if} & r < j \leq 2 r \end{array} \right.$.

The bracket $\{.,.\}$ is a Poisson bracket on $\cO(2r,\nn)$, ie. a bi-derivation which is also a Lie bracket. For all $f \in \cO(2r,\nn)$, let $D_H(f)$ be the derivation $D_H(f): g \mapsto \{ f, g \}$. Then $D_H$ is a morphism of Lie algebras $D_H : \cO(2r,\nn) \to W(2r,\nn)$, the kernel of which is $\FF\, 1$. Let $H(2r,\nn)$ denote its image; then the second derived subalgebra $H(2r,\nn)^{(2)}$ is a simple ideal of dimension $p^{|\nn|} - 2$. Note that $H(2r,\nn)^{(2)}$ is a graded subalgebra of $W(2r,\nn)$; for any $0 \neq \alpha \in \NN^{2r}$, the derivation $D_H(x^{(\alpha)})$ is homogeneous, of degree $|\alpha|-2$.

\subsubsection{} \label{UnEtDemiH}
\begin{thm}
The Lie algebra $H(2r,\nn)^{(2)}$ is not generated by 1,5 elements. \\
\end{thm}

\proof\ We will apply Corollary \ref{corollary1,5} to $L = H(2r,\nn)$ and $H = H(2r,\nn)^{(2)}$. It is easy to check that the centraliser of $H$ in $L$ is trivial. Now let $\Omega \subseteq H(2r,\nn)^{(2)}$ be the open set of elements which have a non-zero component of degree $-1$. We want to show that, for all $D \in \Omega$, there exists $0 \neq \Delta \in H(2r,\nn)_\geq 1$ such that $[D,\Delta] = 0$. 

We choose an element $f \in \cO(2r,\nn)$ such that $D_H(f) = D$. Since $D_H(1) = 0$, we can assume that $f$ has no constant term; also, since the image $D_H(f)$ has a non-zero component of degree $-1$, $f$ must have a non-zero term of degree 1. Then, the lowest degree term of $f^3$ has degree 3, and clearly $\{ f,f^3 \} = 0$. Let $\Delta = D_H(f^3) \in H(2r,\nn)$. Since $f \in \cO(m,\nn)_{\geq 3}$, we have $\Delta \in H(2r,\nn)_{\geq 1}$; also, we have $[D,\Delta] = D_H(\{ f,f^3\}) = 0$. Finally, the fact that $f^3 \not \in \FF$ ensures that $\Delta \neq 0$: this is what we needed to apply Corollary \ref{corollary1,5}, proving the theorem.

\subsection{Lie algebras of type $K$.} $ $\\

In this section, we let $m = 2r + 1 \geq 3$ be an odd integer and $\nn \in \NN^{2r+1}$ a multi-index with non-zero entries. Again, throughout this section the general notations for divided power algebra and special derivations are the same as in \ref{SubSectionW}. 

\subsubsection{}
We recall the constructions from \cite[Section 4.2]{Strade}. As in the case of the Lie algebras of type $H$, we will use the notations $j' = \left\{ \begin{array}{rcl} j+r & \mbox{if} & 1 \leq j \leq r \\ j-r & \mbox{if} & r < j \leq 2 r \end{array} \right.$ and $\sigma(j) = \left\{ \begin{array}{rcl} +1 & \mbox{if} & 1 \leq j \leq r \\ -1 & \mbox{if} & r < j \leq 2 r \end{array} \right.$. (The values of $j'$ and $\sigma(j)$ are not defined for $j = 2r+1$.) 

Define a linear map $D_K: \cO(2r+1,\nn) \to W(2r+1,\nn)$ by the formula:
\begin{equation}
D_K(f) = \sum_{j = 1}^{2r} \Bigg( \sigma(j)\, \d_{j}(f) + x_{j'} \d_m(f) \Bigg) \d_{j'}
+ \left( 2f - \sum_{j=1}^{2r}x_j \d_j(f) \right) \d_{2r+1}.
\end{equation}
For $f,g \in \cO(m,\nn)$, we define the {\it contact bracket of $f$ and $g$} by:
\begin{equation} \label{fundamentalK}
\langle f,g \rangle = D_K(f) (g) - 2\, g\, \d_m(f) \in \cO(m,\nn).
\end{equation}
Then $D_K$ is an injective mapping, and for all $f,g \in \cO(2r+1,\nn)$ we have $[D_K(f) , D_K(g)] = D_K(\langle f,g \rangle)$. The Lie algebra $K(2r+1,\nn)$ is defined to be the image of $D_K$; by construction, the contact bracket $\langle.,.\rangle$ endows $\cO(2r+1,\nn)$ with a Lie algebra structure, and the map $D_K:\cO(2r+1) \to K(2r+1,\nn)$ is an isomorphism of Lie algebras. \\

{\nbf Remark.} Unlike the Poisson bracket on $\cO(2r,\nn)$, the contact bracket is not a biderivation of $\cO(2r+1,\nn)$: for example, we have $\langle1,f \rangle = 2 \d_m(f)$ for all $f \in \cO(2r+1,\nn)$, so that $1$ is not central for the contact bracket.

\subsubsection{} \label{Kdegree}
The derived subalgebra $K(2r+1,\nn)^{(1)}$ is a simple Lie algebra. One has $K(2r+1,\nn)^{(1)} = K(2r+1,\nn)$ unless $m + 3 \equiv 0 \mod p$, in which case $K(2r+1,\nn)^{(1)}$ has codimension 1 in $K(2r+1,\nn)$.

The contact algebra admits a grading such that $K(2r+1,\nn) = K_{-2} \oplus K_{-1} \oplus K_0 \oplus \ldots \oplus K_s$. Unlike the other types, there exist non-zero homogeneous elements of degree -2; in particular, the embedding $D_K:K(2r+1,\nn) \hookrightarrow W(2r+1,\nn)$ does not respect the natural degrees. To emphasise the difference, we will call {\it $K$-degree} the degree defined on $K(2r+1,\nn)$.

Let us describe the $K$-degree in $K(2r+1,\nn)$. For all multi-indices $\alpha \in \NN^m$, we let $\| \alpha \| = \alpha_1 + \ldots +\alpha_{2r} + 2 \alpha_{2r+1}$; this amounts to assigning degree 1 to each variable $x_j$ with $j \leq 2r$, and degree 2 to $x_{2r+1}$. Then, the derivation $D_K(x^{(\alpha)})$ is homogeneous, of $K$-degree $\| \alpha \| - 2$. In particular, if we identify $\cO(2r+1,\nn) \equiv K(2r+1,\nn)$ via $D_K$, we have the following equalities:
\begin{equation} \label{sousespacesK}
K_{-2} = \FF \quad ; \quad K_{-1} = \sum_{j = 1}^{2r} \FF\, x_j \quad ; \quad K_0 = \FF\, x_{2r+1} + \sum_{i,j = 1}^{2r}\, \FF\, x_i\, x_j.
\end{equation}

Note that the $K$-degree is not quite additive for the associative product in $\cO(2r+1,\nn)$: actually, we have $K_i\, K_j \subseteq K_{i + j + 2}$ for all $i,j$.

\subsubsection{}
We will often identify the spaces $\cO(2r+1,\nn)$ and $K(2r+1,\nn)$ by means of the map $D_K$. Each element $f \in \cO(2r+1,\nn)$ induces two linear maps on $\cO(2r+1,\nn)$. First, we have $D_K(f) \in \End\, \cO(2r+1,\nn)$, which is actually a special derivation of the divided power algebra $\cO(2r+1,\nn)$. Second, $\ \ad_K(f)$, the adjoint action of $f$ on the Lie algebra $\cO(2r+1,\nn)$ endowed with the contact bracket: for all $f,g \in \cO(2r+1,\nn)$, we have $\ \ad_K(f)(g) = \langle f,g \rangle$. The link betweeen these mappings is given by the following lemma. \\

\begin{lemma}
Let $f \in \cO(2r+1,\nn)$ be arbitrary. Let $\mu_f \in \End\, \cO(2r+1,\nn)$ be the linear map defined by (associative) multiplication by $f$. Then,
\begin{equation}
\ad_K(f) \circ \mu_f = \mu_f \circ D_K(f) \in \End\, \cO(2r+1,\nn).
\end{equation}
\end{lemma}

\proof\ We use the identity (\ref{fundamentalK}): for all $g \in \cO(2r+1,\nn)$, we have:
\begin{eqnarray*}
\langle f,f g \rangle & = & D_K(f)(fg) - 2 \d_m(f)\, fg \\
 & = & f\, D_K(f)(g) + D_K(f)(f)\, g\, -2 \d_m(f)\, f g \\
  & = &f\, D_K(f)(g) + \langle f,f \rangle\, g =  f\, D_K(f)(g),
\end{eqnarray*}
which means exactly that $ \ad_K(f) \circ \mu_f (g) =  \mu_f \circ D_K(f) (g)$.

\subsubsection{} \label{GenericInvariantsK}
\begin{cor}
Assume that $f$ is invertible in $\cO(2r+1,\nn)$, ie. $f$ has a non-zero constant term. Then $D_K(f)$ and $\ad_K(f)$ are conjugate as linear endomorphisms of $\cO(2r+1,\nn)$. \\
\end{cor}

{\nbf Remark.} Since the conjugating map is given by multiplication by $f$, in particular the invariant subspaces are linked by the relation:
\begin{equation}
\cO(2r+1,\nn)^{\ad_K(f)} = f\, \cO(2r+1,\nn)^{D_K(f)}.
\end{equation}

\subsubsection{} \label{centralisateursK}
\begin{lemma}
Let $D \in K(2r+1,\nn)$ and $K(2r+1,\nn)^D$ be its centraliser.
\begin{enumerate}
\item We have $\dim\, K(2r+1,\nn)^D \geq \min\{ 2r+1 , p \}$.
\item If $D$ has a non-zero term of degree $-2$, then $\dim\, K(2r+1,\nn)_{\geq 1}^D \geq 2$. \\
\end{enumerate}
\end{lemma}

\proof\ We first prove \nr 1. Let $\chi_K(t) = \FF[t]$ be the characteristic polynomial of $D$ acting on $K(2r+1,\nn)$ by the adjoint action. Using \ref{GenericInvariantsK} we can see that it is the same as the characteristic polynomial of $D$ acting as a derivation of $\cO(m,\nn)$, so that by \cite[Theorem 1, (ii)]{Premet_Wn}, it is a $p$-polynomial of degree $p^{|\nn|}$. Thus,
\begin{equation}
\chi_K(t) = f(t)^{p^{|\nn|-d}},
\end{equation}
where $f(t)$ is a separable $p$-polynomial of degree $p^d$ for some $d \in \{ 0,\ldots,|\nn|\}$. Let $\ell \in \{ 0,\ldots , |\nn| - d\}$ be such that the minimal annihilating $p$-polynomial of $D$ is $f_\Lambda(t)^{p^\ell}$.

Let $K(2r+1,\nn)^{(D)}$ be the nilspace of $\ad_K(D)$ in $K(2r+1,\nn)$, so that $K(2r+1,\nn)^D \subseteq K(2r+1,\nn)^{(D)}$. We have $\dim\, K(2r+1,\nn)^{(D)} = p^{|\nn| - d}$, and $D$ induces a nilpotent map of index $\leq p^\ell$ on $K(2r+1,\nn)^{(D)}$. This nilpotent map has the following Jordan decomposition:
\begin{equation}
\left[ \begin{array}{ccc}
 J_1 & & \\
 & \ddots & \\
 & & J_k
\end{array} \right],
\end{equation}
and we have the following information:
\begin{itemize}
\item the sizes of the Jordan cells are $p^\ell \geq |J_1| \geq \ldots \geq |J_k|$,
\item the size of the large matrix is $|J_1| + \ldots + |J_k| = p^{|\nn|-d}$,
\item the number of Jordan cells is $k = \dim\, K(2r+1,\nn)^D$.
\end{itemize}

Assume first that $\ell < |\nn| - d$. By the above properties, we have:
\begin{equation}
p^{|\nn|-d} \leq k |J_1| \leq k p^\ell,
\end{equation}
 so that $\dim\, K(2r+1,\nn)^D = k \geq p^{|\nn|-d-\ell} \geq p$.

Now assume that $\ell = |\nn| - d$. Then the minimal annihilating $p$-polynomial of $\ad_K(D)$ has degree $p^{d+\ell} = p^{|\nn|}$, so that the subspace $\mathfrak{D} = \sum_{j \geq 0} \FF\, \ad_K(D)^{p^j} \subseteq \Der \, K(2r+1,\nn)^{(1)}$ has dimension $|\nn|$. By \cite[Theorem 7.1.2]{Strade}, $K(2r+1,\nn) \subseteq \Der\, K(2r+1,\nn)^{(1)}$ is a subspace of codimension $|\nn| - (2r+1)$, so that $\dim\, K(2r+1,\nn) \cap \mathfrak{D} \geq |\nn| - (|\nn|- 2r-1) = 2r+1$. Since $K(2r+1,\nn) \cap \mathfrak{D} \subseteq K(2r+1,\nn)^D$, we obtain $\dim\, K(2r+1,\nn)^D \geq 2r+1$. \\

We turn to the proof of \nr 2. We identify $K(2r+1,\nn)$ with $\cO(2r+1,\nn)$, endowed with the contact bracket: in this identification we write $D = f \in \cO(2r+1,\nn)$. We assume $D$ has a non-zero term of degree -2, so that equivalently $f$ has a non-zero constant term. By Remark \ref{GenericInvariantsK}, we have $K(2r+1,\nn)^{\ad_K(f)} = f\, K(2r+1,\nn)^{D_K(f)}$. We can check that, more generally:
\begin{equation}
K(2r+1,\nn)^{\ad_K(f)}_{\geq d} = f\, K(2r+1,\nn)^{D_K(f)}_{\geq d},
\end{equation}
for all $d \geq -2$. Hence, to prove \nr 2 we only have to find some non-zero $g \in K(2r+1,\nn)^{D_K(f)}_{\geq 1}$. We insist that we are considering the $K$-grading, which starts in degree $-2$, and not the ordinary grading.

By \nr 1, we have $\dim\,  K(2r+1,\nn)^{D_K(f)} = \dim\, K(2r+1,\nn)^{\ad_K(f)} \geq \min\{ 2r+1,p \} \geq 3$. Therefore, $\dim\, K(2r+1,\nn)^{D_K(f)}_{\geq -1} \geq 2$. 

Let $g \in K(2r+1,\nn)^{D_K(f)}_{\geq -1}$ be a non-zero element. If $g \in K(2r+1,\nn)_{\geq 1}$ we are done. Otherwise, $g$ has a non-zero homogeneous component of degree $-1$ or $0$. Because $D_K(f)$ is an associative derivation, the subring $\FF[g] \subseteq K(2r+1,\nn)^{D_K(f)}$. Using the description of $K_{-1} + K_0$ derived from (\ref{sousespacesK}) we can check that $\dim\, \FF[g] \cap K(2r+1,\nn)_{\geq 1} \geq 2$, and Property \nr 2 follows.

\subsubsection{} \label{UnEtDemiK}
\begin{thm}
The simple Lie algebra $K(2r+1,\nn)^{(1)}$ is not generated by 1,5 elements. \\
\end{thm}

\proof\ We apply Corollary \ref{corollary1,5} to $L = K(2r+1,\nn)$ and $H = K(2r+1,\nn)^{(1)}$. It is quite straightforward to see that the centraliser of $L^H$ is trivial. Then, let $\Omega \subseteq K(2r+1,\nn)^{(1)}$ be the non-empty open subset of derivations with non-zero term of degree $-2$: by Lemma \ref{centralisateursK}, we obtain that for all $D \in \Omega$, the subspace $K(2r+1,\nn)^D \cap K(2r+1,\nn)_{\geq 1} \neq (0)$, which is what we needed to show. The theorem is proved.

\section*{Acknowledgements}

This research was conducted during the author's stay at the Max-Planck Institute in Bonn in January and February 2008. The paper was subsequently written up while visiting the Collaborative Research Centre 701 at the University of Bielefeld in Spring 2008. He would like to express his gratitude towards both institutions for warm hospitality and support.

\end{document}